\input amstex
\magnification=\magstep1
\documentstyle{amsppt}
\advance\hoffset by 0.1 truein
\TagsOnLeft
\font\headlinebf=cmbx12 

\font\smallrm=cmr8
\font\smallit=cmsl8
\nopagenumbers
\document

\

\vskip1.0cm

\subheading\nofrills{\headlinebf Very ampleness for Theta 
on the compactified Jacobian}

\vskip0.5cm

\subheading\nofrills{\bf Eduardo Esteves}\footnote {Supported
by a Starr fellowship from
the MIT Japan Program.}

\subheading{\smallrm Instituto de Matem\'atica Pura e Aplicada, IMPA, Estrada 
Dona Castorina 110, 22460-320 Rio de Janeiro RJ, Brazil}

\vskip0.8cm

\baselineskip=12pt

\subheading\nofrills{\bf 1 Introduction}

\vskip0.4cm

\parindent=0pt

If $X$ is a smooth, complete, connected 
curve over an algebraically 
closed field, then the Jacobian $J_0$, parametrizing
invertible sheaves on $X$ of Euler characteristic $0$, is 
projective and admits a canonical ample divisor $\Theta$, the Theta
divisor. If $g$ denotes the genus of $X$, then $\Theta$ is the 
scheme-theoretic image of the Abel-Jacobi 
morphism $X^{g-1} @>>> J_0$, given by
$$
(p_1,\dots,p_{g-1})\mapsto \Cal O_X(p_1+\dots+p_{g-1}).
$$
It follows from \cite{14, Section 17, p. 163} 
that $3\Theta$ is very ample.

\parindent=16pt

In the singular case, D'Souza has constructed a 
natural compactification $\bar J_0$ for the Jacobian $J_0$ 
of a complete, integral curve over an
algebraically closed field \cite{5}. The scheme $\bar J_0$
parametrizes torsion-free, rank 1 sheaves of Euler
characteristic $0$ on $X$. A natural question 
in this context is whether there is a canonical Cartier 
divisor on $\bar J_0$ extending the notion of the classical Theta
divisor.

The above question was partially and independently
answered in \cite{6} and \cite{19}. In these two works
the same canonical line bundle $\Cal L$
on $\bar J_0$ and the same global section $\theta$ of $\Cal
L$ are defined. For smooth
curves, the zero scheme of $\theta$ is the classical Theta divisor
$\Theta$. In \cite{19} Soucaris shows that the zero
scheme of the restriction of $\theta$ to the maximum
reduced subscheme of $\bar J_0$
is a Cartier divisor. Both \cite{6} and \cite{19} show
that $\Cal L$ is ample. It remains to
determine whether the zero scheme of $\theta$ on $\bar J_0$
is a Cartier divisor in general, and what is the minimum $n$ such that 
$\Cal L^{\otimes n}$ is very ample. 

In this article our main concern is with the latter
question. We will show that $\Cal L^{\otimes n}$ is very ample for $n$ at 
least equal to a specified lower bound (Theorem 7.) 
If $X$ has at most ordinary nodes or cusps as singularities, then 
our lower bound is 3. Our main tool is to use theta sections 
$\theta_E$ associated to vector bundles $E$ on $X$. The
theta sections were used by Faltings \cite{9} to
construct the moduli of semistable vector bundles on a
smooth, complete curve without using Geometric Invariant
Theory (see also \cite{18}.) In a forthcoming work \cite{7}, \cite{8} 
we will apply such
method to construct the compactified Jacobian for families
of reduced curves.

The importance of Theorem 7 is that we obtain a canonical projective 
embedding of $\bar J_0$ in 
$\text{\bf P}(H^0(\bar J_0,\Cal L^{\otimes n}))$, for $n$ minimum 
such that $\Cal L^{\otimes n}$ is very ample. By 
studying the structure of the homogeneous coordinate ring of $\bar J_0$ 
in $\text{\bf P}(H^0(\bar J_0,\Cal L^{\otimes n}))$, 
maybe in a way analogous to Mumford's in \cite{15} and 
\cite{16}, 
we might be able to understand better the algebraic structure of 
$\bar J_0$.

\medpagebreak

\parindent=0pt

\sl Notation. \rm We will often deal with parameter spaces, that is, spaces
whose points are classes representing certain objects. In
such context, we will employ the usual bracket notation
$[F]$ for the point representing the object $F$. 
If $E$ is a vector bundle on a scheme $Y$, we denote by 
$\text{\bf P}_Y(E)$ the corresponding projective bundle over $Y$. 
By a point we mean a closed point.

\vskip0.8cm

\subheading\nofrills{\bf 2 The compactified Jacobian}

\vskip0.4cm

\headline{\smallrm 2 \  \  \  \  \  \  \  \  \  \  \  
\  \  \  \  \  \  \  \  \  \  \  \  \  \  \  \  \  \  \  \  \  \  \  \  \  
\  \  \  \  \  \  \  \  \  
\  \  \  \  \  \  \  \  \  
\  \  \  \  \  \  \  \  \  \  \  \  \  \  E. Esteves}

Let $X$ be a complete, integral 
curve over an algebraically closed field $k$. Denote by 
$g$ the arithmetic genus of $X$, and by $\omega$ the 
dualizing sheaf on $X$. A coherent sheaf $I$ on $X$ is \sl torsion-free 
\rm if $I_x$ is a torsion-free $\Cal O_x$-module for every $x\in X$. 
A coherent sheaf $I$ on $X$ is \sl rank 1 \rm if $I$ is generically 
invertible. By \cite{4, p. 96}, the sheaf 
$\omega$ is torsion-free, rank 1. Fix an ample line
bundle $\Cal O_X(1)$ on $X$. For every coherent sheaf
$F$ on $X$, let $\text{Quot}^{p(t)}_X(F)$ denote
Grothendieck's Quot-scheme \cite{10}, parametrizing quotients of $F$
with Hilbert polynomial $p(t)$ with respect to $\Cal
O_X(1)$. We will drop the superscript $p(t)$ whenever it
is not important. 

\parindent=16pt

For every integer $d$, let $\bar {\text{\bf J}}_d$ 
denote the \sl compactified Jacobian 
functor. \rm For each $k$-scheme $S$, the set $\bar {\text{\bf J}}_d(S)$ 
consists of equivalence classes of $S$-flat coherent sheaves $\Cal I$ on 
$X\times S$ such that $\Cal I(s)$ is torsion-free, rank 1 
of Euler characteristic $d$ for every $s\in S$. (Two sheaves $\Cal I_1$ and 
$\Cal I_2$ are called equivalent if there is an invertible sheaf $N$ on $S$ 
such that $\Cal I_1\cong\Cal I_2\otimes N$.) 
D'Souza \cite{5} and Altman and Kleiman \cite{3}, \cite{4} have shown that 
$\bar {\text{\bf J}}_d$ is represented by a (projective) 
scheme $\bar J_d$, the 
\sl compactified Jacobian. \rm Here we present yet another 
proof of the representability of $\bar {\text{\bf J}}_d$ by a scheme, a 
proof more suitable for treating the question of very ampleness in 
Section 4.

For every torsion-free, rank 1 sheaf $I$ on $X$, let
$$
e(I):=\max_{x\in X} \dim_k I(x).
$$
Since $I$ is generically invertible, then $1\leq e(I)<\infty$.

\proclaim\nofrills{\bf Proposition 1} \  Let $I$ be 
a torsion-free, rank 1 sheaf on $X$ of
Euler characteristic $d$. Then, for every integer 
$r\geq \max(e(\underline{Hom}
_X(I,\omega)),2)$, there is a vector bundle
$E$ on $X$ of rank $r$ and degree $-rd-1$ such that:
\smallpagebreak

\parindent=0pt

(i) \  $h^0(X,I\otimes E)=0$ and $h^1(X,I\otimes E)=1$;

(ii) the unique (modulo $k^*$)
non-zero homomorphism $I @>>>
E^*\otimes\omega$ is an embedding with torsion-free cokernel.
\endproclaim

\parindent=16pt

\demo{Proof} Let $m>>0$ be an integer such that 
$H^0(X,I(-m))=0$ and 
$\underline{Hom}_X(I,\omega)(m)$ is generated by global sections. 
Since $r\geq \max(e(\underline{Hom}_X(I,\omega)),2)$ and $k$ is infinite, 
then there is a surjection $p\: \Cal O_X^{\oplus r} \twoheadrightarrow
\underline{Hom}_X(I,\omega)(m)$. Applying $\underline{Hom}_X(\cdot,\omega)$ 
to $p$, we obtain an embedding $I(-m) \hookrightarrow \omega^{\oplus r}$, 
whose cokernel is torsion-free since $\underline{Ext}^1_X(F,\omega)=0$ for 
every torsion-free sheaf $F$ on $X$ \cite{4, p. 96}. 
Twisting by $\Cal O_X(m)$ and letting $E:=\Cal O_X(-m)^{\oplus r}$, we 
get that $H^0(X,I\otimes E)=0$ and there is a short exact sequence on $X$ 
of the form 
$$
0 @>>> I @>\mu>> E^*\otimes\omega @>q>> C @>>> 0,
$$
where $C$ is torsion-free. 

Let $h:=h^1(X,I\otimes E)$. 
If $h=1$, then the proposition is proved. We will show by
descending induction on $h$ that we can choose $E$ as in the 
above paragraph with $h=1$. Suppose $h>1$. 
Let $\lambda\: I @>>> E^*\otimes\omega$ be a 
homomorphism that is not a multiple of $\mu$. Since $I$ is simple by 
\cite{4, Lemma 5.4, p. 83}, 
then the composition $\rho:=q\circ\lambda$ is not zero. 
Since $C$ is torsion-free, then there is a regular point $x\in X$
such that $\rho(x)\neq 0$. Let $\sigma\: C(x) @>>> I(x)$ be a splitting 
for $\rho(x)$. Let
$$
F:=(\ker (E^* @>>> E^*(x) @>q(x)>> C(x) @>\sigma >> I(x)))^*.
$$
(We implicitly chose a trivialization of $\omega$ at $x$. Any other choice of 
trivialization yields the same subsheaf $F$.) 
Then $F$ is a vector bundle of $\deg F=\deg E+1$ and rank $r$. 
By definition of $F$, we have that $\mu$ factors
through an embedding $\mu'\: I @>>> F^*\otimes\omega$, but
$\lambda$ does not. Thus $h^1(X,I\otimes F)<
h^1(X,I\otimes E)$. Since $\deg F=\deg E+1$, then
$H^0(X,I\otimes F)=0$. It is clear that the
cokernel of $\mu'$ is torsion-free. The induction proof is complete.
\qed
\enddemo

\headline{\smallrm Very ampleness for Theta on the compactified Jacobian 
\  \  \  \  \  \  \  \  \  \  \  \  \  \  \  \  \  \  \  \  
\  \  \  \  \  \  \  \  \  \  \  \  \  \  \  \  \  \  \  \  \  \  \  \  \  
\  \  \  \  \  \  \  \  \  \  \  \  \  \  \  \  3}

\proclaim\nofrills{\bf Corollary 2} \  
The functor $\bar {\text{\bf J}}_d$ is representable by a scheme.
\endproclaim

\demo{Proof} First note that properties (i) and (ii) in the statement of 
Proposition 1 are open on $I$. 
More precisely, given a vector bundle $E$ on $X$ of rank $r$ 
and degree $-rd-1$, the subfunctor $\text{\bf U}_E\subseteq \bar {\text{\bf 
J}}_d$, parametrizing sheaves $I$ satisfying properties (i) and (ii) in the
statement of Proposition 1, is open. By Proposition 1, the subfunctors 
$\text{\bf U}_E$ cover $\bar {\text{\bf J}}_d$. Thus to show that 
$\bar {\text{\bf J}}_d$ is representable we need only show that each 
$\text{\bf U}_E$ is representable.

Fix a vector bundle $E$ on $X$ of rank $r$ and degree $-rd-1$. Let
$$
V\subseteq \text{Quot}_X(E^*\otimes\omega)
$$
be the open subscheme parametrizing those quotients 
$q\: E^*\otimes\omega @>>>
G$ such that both $G$ and $\ker(q)$ are torsion-free, 
$\ker(q)$ has rank 1, 
$$
h^0(X,\  \ker(q)\otimes E)=0 \text{\  \  and \  \  } 
h^1(X,\  \ker(q)\otimes E)=1.\tag{2.1}
$$
There is a morphism of functors $V @>>> \text{\bf U}_E$ 
sending a quotient $[q]\in V$ to its kernel, 
$[\ker(q)]\in \text{\bf U}_E$. 
It follows from (2.1) that the latter morphism is an isomorphism. 
The proof is complete.\qed 
\enddemo

We will say that a vector bundle $E$
on $X$ of rank $r$ and degree $-rd-1$ satisfying properties (i)
and (ii) in the statement of Proposition 1 \sl represents $I$.
\rm We remark that the property of representing $I$ is open.

\vskip0.8cm

\subheading\nofrills{\bf 3 The Theta divisor}

\vskip0.4cm

\parindent=0pt

Assume from now on that $g>0$. 
Let $\Cal I$ be a universal relatively torsion-free, rank 1 sheaf on
$X\times \bar J_0$ over $\bar J_0$. Denote by 
$p\: X\times\bar J_0 @>>> \bar
J_0$ the projection map. Define
$$
\Cal L:=(\det Rp_*(\Cal I))^{-1},
$$
where $\det Rp_*$ denotes the determinant of cohomology
associated with the projection $p$. (For a brief description of
$Rp_*$, see \cite{6} or \cite{19}; for a
more in-depth development of the theory of determinants,
see \cite{12}.) Since the sheaf $\Cal I$ has relative
Euler characteristic $0$ over $\bar J_0$, then $\Cal L$ is
independent on the choice of a universal sheaf $\Cal I$, and there is
a canonical global section $\theta$ of $\Cal L$
whose zero scheme $\Theta$ parametrizes torsion-free,
rank 1 sheaves $I$ of Euler characteristic $0$ on $X$ such that
$$
h^0(X,I)=h^1(X,I)\neq 0.
$$
Equivalently, by Serre's duality, $\Theta$ 
consists of the torsion-free, rank 1 sheaves of Euler
characteristic $0$ that can be embedded into the 
dualizing sheaf $\omega$. In other words, $\Theta$ is (set-theoretically) 
the image of the 
$(g-1)$-th component of the Abel-Jacobi map:
$$
\Cal A^{g-1}\:\text{Quot}^{g-1}_X(\omega) @>>> \bar J_0,
$$
where $\Cal A^{g-1}$ sends a quotient
$[q]\in\text{Quot}^{g-1}_X(\omega)$ 
to its kernel, $[\ker(q)]\in \bar J_0$ (cf.
\cite{4, p. 87}.) We say that $\Cal L$ is the \sl Theta line 
bundle\rm , and $\Theta$ is 
the \sl Theta divisor \rm (even though it is not known 
whether $\Theta$ is actually a Cartier divisor in general.)

\headline{\smallrm 4 \  \  \  \  \  \  \  \  \  \  \  
\  \  \  \  \  \  \  \  \  \  \  \  \  \  \  \  \  \  \  \  \  \  \  \  \  
\  \  \  \  \  \  \  \  \  \  \  \  \  \  \  \  \  \  \  \  \  \  \  \  \  
\  \  \  \  \  \  \  \  E. Esteves}

\parindent=16pt

If $X$ is smooth, then 
$$
\text{Quot}^{g-1}_X(\omega)\cong\text{Hilb}^{g-1}_X=\text{Symm}^{g-1}(X),
$$
where $\text{Hilb}^{g-1}_X:=\text{Quot}^{g-1}_X(\Cal O_X)$
is the Hilbert scheme, parametrizing $(g-1)$-uples of points
in $X$, and $\text{Symm}^{g-1}(X)$ is the symmetric
product of $(g-1)$ copies of $X$. Hence $\Theta$
corresponds to the classical Theta divisor (cf. Section 1.)

Assume that $X$ is locally planar, that
is, that the embedding dimension of each point of $X$ is at
most 2. Equivalently, assume that $X$ can be embedded into a quasi-projective
smooth surface \cite{2}. Then $\text{Quot}^{g-1}_X(\omega)$ 
and $\bar J_0$ are integral, 
local complete intersections of dimensions $g-1$ and $g$, respectively 
(Since locally planar curves
are Gorenstein, then
$\text{Quot}^{g-1}_X(\omega)\cong\text{Hilb}^{g-1}_X$, and
thus our statement follows from \cite{1, Cor. 7 and Thm. 9}.) 
In this case, $\Theta$ is an irreducible, local complete
intersection, effective Cartier divisor 
on $\bar J_0$. Moreover, 
it is clear that $\Cal A^{g-1}$ is an isomorphism over the open 
subscheme of $\Theta$ parametrizing torsion-free, rank 1 sheaves $I$ with 
$h^1(X,I)=1$. From \cite{4, Prop. 3.5.ii, p. 76}, 
this open subscheme is dense. Since $\Theta$ is
Cohen-Macaulay 
and irreducible, and $\text{Quot}^{g-1}_X(\omega)$ is
integral, then $\Theta$ is also 
integral. We observe that the assumption that $X$ is
locally planar is essential in the above argument.
If $X$ is not locally planar, then $\bar J_0$ is not
irreducible (cf. \cite{11} or \cite{17, Thm. A}), 
and may have dimension greater
than $g$ (cf. \cite{1, Ex. 13, p. 10}.)

We observe that the above notions and arguments can
be extended to families of integral, complete curves 
without difficulty \cite{6}. Moreover, the formation of $\theta$
and $\Cal L$ commutes with base change, since so does
the determinant of cohomology. From this
observation it follows that Poincar\'e's
formula holds for locally planar curves. Namely, we claim
that, if $X$ is locally planar, then the 
self-intersection $\Theta^g$ is equal to $g!$. 
In fact, the claim is known for smooth curves 
\cite{13, \S 2}. Since every locally planar curve is part of a
family whose general member is a smooth curve, then we may
apply the principle of conservation of intersection number
to prove our claim.

As we have already remarked, it is not known whether 
$\Theta$ is always a Cartier divisor. 
Nevertheless, Soucaris showed that the zero scheme
of the restriction of the canonical section $\theta$ to
the maximum reduced subscheme of $\bar J_0$ is a Cartier
divisor \cite{19, Thm. 8, p. 236}.

\vskip0.8cm

\headline{\smallrm Very ampleness for Theta on the compactified Jacobian 
\  \  \  \  \  \  \  \  \  \  \  \  \  \  \  \  \  \  \  \  
\  \  \  \  \  \  \  \  \  \  \  \  \  \  \  \  \  \  \  \  \  \  \  \  \  
\  \  \  \  \  \  \  \  \  \  \  \  \  \  \  \  5}

\subheading\nofrills{\bf 4 Very ampleness}

\vskip0.4cm

\parindent=0pt

Recall the notations of Section 3. If $E$ is a vector bundle
on $X$ with $\deg E=0$, then $\Cal I\otimes E$
has relative Euler characteristic $0$ over $\bar J_0$.
Therefore, the invertible sheaf
$$
\Cal L_E:=(\det Rp_*(\Cal I\otimes E))^{-1}
$$
on $\bar J_0$ has a canonical global section $\theta_E$,
whose zero scheme $\Theta_E$ parametrizes torsion-free,
rank 1 sheaves $I$ of Euler characteristic $0$ on $X$ such
that
$$
h^0(X,I\otimes E)=h^1(X,I\otimes E)\neq 0.
$$
As before, $\Cal L_E$ and $\theta_E$ are
independent on the choice of a universal sheaf $\Cal I$.

\parindent=16pt

\proclaim\nofrills{\bf Lemma 3} \  Let $E$ and $F$ be vector bundles on $X$
of same rank and degree $0$. If $\det E\cong\det F$, then
$\Cal L_E\cong\Cal L_F$.
\endproclaim

\demo{Proof} By Seshadri in \cite{18, Lemma 2.5, p.
165}.\qed
\enddemo

By Lemma 3, if $E$ is a vector bundle on $X$ of rank $n$
and $\det E\cong\Cal O_X$, then $\Cal L_E\cong\Cal
L^{\otimes n}$. Thus we may 
consider $\theta_E$ as a global section of $\Cal
L^{\otimes n}$ under the latter isomorphism. 
In this case we say that $\theta_E$ is a \sl theta section 
of degree $n$. \rm We now have a convenient way 
to produce sections of powers of
$\Cal L$.

For every integer $d$, let $J^d$ be the Jacobian of $X$,
parametrizing invertible sheaves of degree $d$ on $X$.
Recall that $J^d$ is connected, quasi-projective and smooth. 

\proclaim\nofrills{\bf Lemma 4} \  Let $n\geq 2$. For each $i=1,\dots, n$,
let $d_i$ be an integer and 
$U_i\subseteq J^{d_i}$ be a non-empty, open subset.
Let $L$ be an invertible sheaf of degree $d_1+\dots+d_n$. Then there
are points $[L_i]\in U_i$ for every $i=1,\dots, n$
such that
$$
L\cong \bigotimes_{i=1}^n L_i.
$$
\endproclaim

\demo{Proof} Consider the morphism $\phi\:
U_1\times\dots\times U_{n-1} @>>> J^{d_n}$, given by
$$
([M_1],\dots,[M_{n-1}]) \mapsto [L\otimes M_1^{-1}\otimes\dots\otimes
M_{n-1}^{-1}].
$$
It is clear that the image $V$ of $\phi$ is
open in $J^{d_n}$. Since $J^{d_n}$ is irreducible, then $V\cap
U_n\neq\emptyset$. Thus there is a point $[L_i]\in U_i$
for each $i=1,\dots, n$ such that 
$$
L_n\cong L\otimes L_1^{-1}\otimes\dots L_{n-1}^{-1}.
$$
The proof is complete.
\qed
\enddemo

\proclaim\nofrills{\bf Theorem 5} \  The sheaf $\Cal L^{\otimes n}$ is
generated by global sections if $n\geq 2$.
\endproclaim

\demo{Proof} Fix $n\geq 2$. Let $I$ be a torsion-free,
rank 1 sheaf on $X$ of Euler characteristic 0. We will show that
there is a vector bundle $E$ on $X$ of rank $n$ and $\det
E\cong\Cal O_X$ such that
$$
h^0(X,I\otimes E)=h^1(X,I\otimes E)=0.\tag{5.1}
$$
In this case, the section $\theta_E$ generates $\Cal
L^{\otimes n}$ at
$[I]$, thereby proving the theorem.

By the proof of \cite{19, Prop. 7, p. 235}, there is an 
invertible sheaf $L$ on $X$ of degree $0$ such that
$$
h^0(X,I\otimes L)=h^1(X,I\otimes L)=0.
$$
By semicontinuity, there is an open, dense subset
$U\subseteq J^0$, containing $[L]$, such that if $[M]\in U$, then
$$
h^0(X,I\otimes M)=h^1(X,I\otimes M)=0.
$$
From Lemma 4, with $U_i:=U$ for every $i=1,\dots, n$, 
there are invertible sheaves $M_1,\dots, M_n$ of degree
$0$ on $X$ such that
$$
h^0(X,I\otimes M_i)=h^1(X,I\otimes M_i)=0
$$
for every $i=1,\dots, n$, and
$$
M_1\otimes\dots \otimes M_n\cong \Cal O_X.
$$
If we now let $E:=M_1\oplus \dots \oplus M_n$, then $E$
satisfies (5.1) and $\det E\cong\Cal O_X$. 
The proof is complete.\qed
\enddemo

\headline{\smallrm 6 \  \  \  \  \  \  \  \  \  \  \  
\  \  \  \  \  \  \  \  \  \  \  \  \  \  \  \  \  \  \  \  \  \  \  \  \  
\  \  \  \  \  \  \  \  \  \  \  \  \  \  \  \  \  \  \  \  \  \  \  \  \  
\  \  \  \  \  \  \  \  \  \  \  \  \  \  \  E. Esteves}

Soucaris had used \cite{19, Prop. 7, p. 235} to show 
that the pullback of $\Cal
L^{\otimes 2}$ to the normalization of $\bar J_0$ is
generated by global sections \cite{19, Prop. 9, p. 236}.

If $S$ is a
$k$-scheme and $\Cal F$ is a vector bundle on $X\times S$
of relative degree $d$ over $S$, then we denote by
$\pi_{\Cal F}\: S @>>> J^d$ the determinant morphism,
mapping $s\in S$ to $[\det\Cal F(s)]\in J^d$.

\proclaim\nofrills{\bf Lemma 6} \  Let $F_1,\dots, F_n$ be vector bundles
on $X$ of same rank $r$ and same degree $d$. Then 
there are a connected, smooth $k$-scheme $S$ and a vector bundle $\Cal F$ on
$X\times S$ such that $\pi_{\Cal F}$ is smooth, and
$F_i\cong\Cal F(s_i)$ for some $s_i\in S$, for each $i=1,\dots,n$.
\endproclaim

\demo{Proof} Let $m>>0$ be such that $F_i(m)$ is
generated by global sections for every $i=1,\dots,n$. 
Since $k$ is infinite, then there is an exact sequence of
the form
$$
0 @>>> \Cal O_X(-m)^{\oplus r-1} @>>> F_i @>>> (\det F_i)((r-1)m)
@>>> 0\tag{6.1}
$$
for each $i=1,\dots, n$. Let $\Cal P$ be a universal sheaf on
$X\times J^d$. Let $p\: X\times J^d @>>> J^d$ denote the
projection map, and let $\Cal V:=R^1p_*(\Cal
P^{-1}(-rm))^{\oplus r-1}$. Choose $m>>0$ such that $\Cal V$ is
locally free, and let $T:=\text{\bf P}_{J^d}(\Cal V^*)$. 
Since $\Cal V$ is locally free, then $T$ is smooth over
$J^d$. Since $J^d$ is connected, smooth and
quasi-projective, then so is $T$. The
scheme $T$ parametrizes $\Cal O_X$-module extensions of
$L((r-1)m)$ by $\Cal O_X(-m)^{\oplus r-1}$ for invertible sheaves $L$ on $X$ 
of degree $d$. Thus there is $s_i\in T$ 
corresponding to (6.1) for each $i=1,\dots, n$. Since $T$ is
quasi-projective, then there is an affine open subscheme
$S\subseteq T$ 
containing $s_1,\dots s_n$. Since $S$ is affine, then 
$$
\Cal V(S)=\text{Ext}^1_{X\times S}(\left.\Cal P\right|_{X\times S}((r-1)m),
\Cal O_{X\times S}(-m)^{\oplus r-1}).
$$
Let $q\: \Cal V^*_T @>>> \Cal Q$ be the universal quotient
on $T$ over $J^d$. Then $q$ induces an extension of the form
$$
0 @>>> \Cal O_{X\times S}(-m)^{\oplus r-1}\otimes\Cal Q @>>> \Cal F @>>>
\left.\Cal P\right|_{X\times S}((r-1)m) @>>> 0
$$
on $X\times S$ that specializes to (6.1) over $s_i$, for each $i=1,\dots,n$.
By construction,
$\pi_{\Cal F}$ is equal to the restriction to $S$ of
the structure morphism $T @>>> J^d$. Thus $\pi_{\Cal F}$
is smooth. The proof is complete.\qed
\enddemo

\headline{\smallrm Very ampleness for Theta on the compactified Jacobian 
\  \  \  \  \  \  \  \  \  \  \  \  \  \  \  \  \  \  \  \  
\  \  \  \  \  \  \  \  \  \  \  \  \  \  \  \  \  \  \  \  \  \  \  \  \  
\  \  \  \  \  \  \  \  \  \  \  \  \  \  \  \  \  7}

Let $e_X:=\max_I e(I)$, 
where the maximum runs over all torsion-free, rank 1 sheaves on $X$. If 
$S$ is a $k$-scheme, we say that an $S$-flat coherent sheaf $\Cal C$ on 
$X\times S$ is \sl relatively torsion-free \rm if $\Cal C(s)$ is 
torsion-free for every $s\in S$.

\proclaim\nofrills{\bf Theorem 7} \  The sheaf $\Cal L^{\otimes n}$ is
very ample for every $n\geq \max(e_X,2)+1$.
\endproclaim

\demo{Proof} Fix $n\geq \max(e_X,2)+1$. By Theorem 5,
the sheaf $\Cal
L^{\otimes n}$ is generated by global sections. We need
only show that $\Cal L^{\otimes n}$ separates points and
tangent vectors on $\bar J_0$. The former is 
Step 1, while the latter is Step 2 below.

\medpagebreak

\parindent=0pt

\sl Step 1: Let $I_1$ and $I_2$ be non-isomorphic 
torsion-free, rank 1 sheaves on $X$ of Euler
characteristic $0$. Then there is a
vector bundle $E$ on $X$ of rank $n$ and $\det E\cong\Cal O_X$
such that $\theta_E([I_1])\neq 0$, but $\theta_E([I_2])=0$.

\medpagebreak

Proof of Step 1: \rm By Proposition 1, since $n\geq \max(e_X,2)+1$, 
there is a vector bundle
$F_i$ on $X$ of rank $n-1$ and degree $-1$ representing 
$I_i$ for each $i=1,2$. From Lemma 6, since the property of representing a
torsion-free, rank 1 sheaf is open, we may 
assume that there are a non-empty, connected, smooth
$k$-scheme $S$, and a vector bundle $\Cal F$ on $X\times S$ 
of rank $n-1$ and relative degree $-1$ over $S$ such that 
the determinant morphism $\pi_{\Cal F}$ is smooth, and 
$\Cal F(s)$ represents both $I_1$ and $I_2$ for every
$s\in S$.

\parindent=16pt

By replacing $S$ with an open, dense subscheme if necessary, we
may assume that for each $i=1,2$ there is an exact
sequence
$$
0 @>>> I_i\otimes\Cal O_S @>\lambda_i>> \Cal
F^*\otimes\omega @>q_i>> \Cal C_i @>>> 0
$$
on $X\times S$, where $\Cal C_i$ is a relatively torsion-free 
sheaf over $S$. If the composition $\rho:= 
q_2\circ \lambda_1$ were zero over a certain $s\in
S$, then $\lambda_1(s)$ would factor through $I_2$, and
since $\chi(I_1)=\chi(I_2)$ we would have that $I_1\cong
I_2$. Thus $\rho\: I_1\otimes\Cal O_S @>>> \Cal C_2$ is an
embedding with $S$-flat cokernel. Since $\Cal C_2$ is
relatively torsion-free, by replacing $S$ with an
open, dense subscheme if necessary, we may assume that there is
a regular point $x\in X$ such that $\rho(x)\:
I_1(x)\otimes\Cal O_S @>>> \Cal C_2(x)$ is an embedding with free 
cokernel. Let $\sigma\: \Cal C_2(x) @>>> I_1(x)\otimes\Cal O_S$ be a 
splitting for $\rho(x)$. Let 
$$
\Cal G:=(\ker (\Cal F^* @>>> \Cal F^*(x) @>q_2(x)>>
\Cal C_2(x) @>\sigma >> I_1(x)\otimes\Cal O_S))^*.
$$
(As in the proof of Proposition 1, we implicitly chose a trivialization of 
$\omega$ at $x$.) Then $\Cal G$ is a vector bundle on $X\times S$ of rank
$n-1$ and relative degree 0 over $S$. Moreover, 
$\det\Cal G(s)\cong\det \Cal F(s)\otimes\Cal O_X(x)$ for
every $s\in S$. Thus the determinant morphism $\pi_{\Cal
G}$ is also smooth. In addition, $\lambda_2$ factors through $\Cal
G^*\otimes\omega$, but $\lambda_1(s)$ does not factor
through $\Cal G^*(s)\otimes\omega(s)$ for any $s\in S$. Thus
$$
h^0(X,I_1\otimes\Cal G(s))=0, \text{\  \  but \  \  } 
h^0(X,I_2\otimes\Cal G(s))\neq 0
$$
for every $s\in S$.

By the proof of \cite{19, Prop. 7, p. 235} (see the proof of Theorem
5), there is an open dense subset $U\subseteq J^0$ such that 
$$
h^0(X,I_1\otimes L)=h^0(X,I_2\otimes L)=0
$$
for every $[L]\in U$. By Lemma 4 applied to
$U_1:=\pi_{\Cal G}(S)$ and $U_2:=U$, there are $s\in S$ and
$[L]\in U$ such that
$$
(\det \Cal G(s))\otimes L\cong\Cal O_X.
$$
It is clear that $E:=\Cal G(s)\oplus L$ satisfies $\det E\cong\Cal O_X$ and
$$
h^0(X,I_1\otimes E)=0, \text{\  \  but \  \  }
h^0(X,I_2\otimes E)\neq 0.
$$
Thus $\theta_E([I_1])\neq 0$, but $\theta_E([I_2])=0$. 
The proof of Step 1 is complete.\qed

\medpagebreak

\headline{\smallrm 8 \  \  \  \  \  \  \  \  \  \  \  
\  \  \  \  \  \  \  \  \  \  \  \  \  \  \  \  \  \  \  \  \  \  \  \  \  
\  \  \  \  \  \  \  \  \  \  \  \  \  \  \  \  \  \  \  \  \  \  \  \  \  
\  \  \  \  \  \  \  \  \  \  \  \  \  E. Esteves}

\parindent=0pt

\sl Step 2: Let $I$ be a torsion-free, rank 1 sheaf
on $X$ with $\chi(I)=0$. Let $v\in T_{\bar J_0,[I]}$ be a
non-zero tangent vector on $\bar J_0$ at $[I]$. Then 
there is a vector bundle $E$ on $X$ of rank $n$ 
and $\det E\cong\Cal O_X$ such that $\Theta_E$ contains 
$[I]$ but not $v$.

\medpagebreak

Proof of Step 2: \rm As in Step 1, we may assume that
there are a non-empty, connected, smooth $k$-scheme $S$,
and a vector bundle $\Cal F$ on $X\times S$ of rank $n-1$ and
relative degree $-1$ over $S$, such that the determinant
morphism $\pi_{\Cal F}$ is smooth, and $\Cal F(s)$
represents $I$ for every $s\in S$.

\parindent=16pt

By replacing $S$ with an open, dense subscheme if necessary, we may assume
that there is an exact sequence
$$
0 @>>> I\otimes\Cal O_S @>\lambda>> \Cal F^*\otimes\omega
@>q>> \Cal C @>>> 0
$$
on $X\times S$, where $\Cal C$ is a relatively 
torsion-free sheaf over $S$. 
By the proof of Corollary 2, we have natural 
identifications
$$
T_{\bar J_0,[I]}=T_{Q(s),[q(s)]}=Hom_X(I,\Cal C(s))\tag{7.1}
$$
for every $s\in S$, where $Q(s):=\text{Quot}_X(\Cal F^*(s)\otimes\omega)$. 
So there is a homomorphism $\nu\:
I\otimes\Cal O_S @>>> \Cal C$ such that $\nu(s)=v$ under the identification 
(7.1) for every $s\in S$. Since $v\neq 0$, then $\nu$ is an embedding with
$S$-flat cokernel. Since $\Cal C$ is relatively
torsion-free, by replacing $S$ with an open, dense subscheme if
necessary, there is a regular point $x\in X$ such that $\nu(x)$ is an
embedding with free cokernel. Let $\sigma\: \Cal C(x) @>>> 
I(x)\otimes\Cal O_S$ be a splitting for $\nu(x)$. Let
$$
\Cal G:=(\ker (\Cal F^* @>>> \Cal F^*(x) @>q(x)>>
\Cal C(x) @>\sigma >> I(x)\otimes\Cal O_S))^*.
$$
(As in the proof of Step 1, we implicitly chose a trivialization of 
$\omega$ at $x$.) Then $\Cal G$ is a vector bundle on $X\times S$ of rank 
$n-1$ and relative degree 0 over $S$. Moreover, $\det\Cal
G(s)\cong \det\Cal F(s)\otimes\Cal O_X(x)$ for every $s\in
S$. Thus the determinant morphism $\pi_{\Cal G}$ is
smooth. In addition, $\lambda$ factors through $\Cal
G^*\otimes\omega$. Thus $[I]\in\Theta_{\Cal G(s)}$ for
every $s\in S$. On the other hand, 
since $\nu(x)$ is an embedding, then $v$
does not belong to $\Theta_{\Cal G(s)}$
for any $s\in S$.

The reader is invited to repeat the argument in the last
paragraph of the proof of Step 1 to finish the proof of Step 2. The proof 
of Theorem 7 is complete.\qed
\enddemo

\parindent=0pt

\headline{\smallrm Very ampleness for Theta on the compactified Jacobian 
\  \  \  \  \  \  \  \  \  \  \  \  \  \  
\  \  \  \  \  \  \  \  \  \  \  \  \  \  \  \  \  \  \  \  \  \  \  \  \  
\  \  \  \  \  \  \  \  \  \  \  \  \  \  \  \  \  9}

\sl Remark 8. \rm Let $x\in X$. 
Let $\overline{\Cal O}_x$ denote the normalization of 
$\Cal O_x$. Let $\delta_x$ denote the length of $\overline{\Cal O}_x/
\Cal O_x$. If $I$ is a torsion-free, rank 1 module over $\Cal O_x$, 
then it is easy to show that $I$ is isomorphic to a submodule of 
$\overline{\Cal O}_x$ containing $\Cal O_x$. Thus 
$$
\dim_k I(x)\leq \delta_x+1.\tag{8.1}
$$
If the conductor, $\Cal C_x:=(\Cal O_x:
\overline{\Cal O}_x)$, is the maximal ideal $m_x$ of $\Cal O_x$, then 
equality in (8.1) is achieved for $I=\overline{\Cal O}_x$ only; 
otherwise the inequality (8.1) is always strict. Let
$$
\delta_X:=\max_{x\in X} \delta_x.
$$
Since $X$ is generically non-singular, then $\delta_X<\infty$. It follows 
from (8.1) that $e_X\leq\delta_X+1$.

\parindent=16pt

Theorem 7 states that 
$\Cal L^{\otimes 3}$ is very ample if $e_X\leq 2$. This is the case for 
$X$ non-singular, or with at most 
ordinary nodes or cusps as singularities, as $\delta_X\leq 1$. It is clear 
that if $e_X\leq 2$ then $X$ is locally planar. If $\delta_X=2$, 
then $e_X\leq 2$ if and only if $X$ is 
locally planar. If $\delta_X=3$, then $e_X\leq 2$ if and only if 
$X$ is locally planar and $m_x^2\neq\Cal C_x$ for every $x\in X$. Note 
that the planar curve $X\subseteq\text{\bf P}^2_k$, 
given as the zero scheme of $u^3w-v^4$, has $e_X=3$.

\medpagebreak

\parindent=0pt

\sl Question 9. \rm It follows from the proof of the theorem in 
\cite{14, Section 17, p. 163} that, if $X$ is smooth, 
then $3\Theta$ is very ample, and the sections 
$\theta_E$ associated to completely decomposable vector bundles $E$ (that is: 
vector bundles $E$ of the form $L_1\oplus L_2\oplus L_3$, 
where $L_i$ is an invertible sheaf of degree 0 
for $i=1,2,3$, and $L_1\otimes L_2\otimes L_3\cong\Cal O_X$), 
are enough to embed $J_0$ into a projective space. We might ask: for which 
integral curves $X$ are such 
sections enough to embed $\bar J_0$ into a projective space? 
The proof of Theorem 7 shows that the sections $\theta_E$ associated to 
vector bundles $E$ of the form $F\oplus L$, 
where $F$ is a vector bundle of rank $\max(e_X,2)$ and 
degree 0, the sheaf $L$ is invertible of degree 0 and 
$(\det F)\otimes L\cong\Cal O_X$, 
are enough to embed $\bar J_0$ into a projective 
space.

\medpagebreak

\sl Example 10. \rm Let $X$ be a complete, integral curve of arithmetic genus
$g=1$. As a subset, $\Theta$ is the
locus of torsion-free, rank 1 sheaves $I$ with Euler
characteristic 0 such that $h^0(X,I)>0$. Since $\chi(\Cal
O_X)=0$, then any non-zero section $\Cal O_X @>>> I$ must
be an isomorphism. Since $\Theta$ is integral by
Subsection 3, then $\Theta=[\Cal O_X]$, as Cartier
divisors of $\bar J_0$. 

\parindent=16pt

By \cite{3, Ex. 8.9.iii, p. 109}, the first component of the
Abel-Jacobi map,
$$
\aligned 
\Cal A^1\: X @>>> & \bar J_{-1}\\
x\   \mapsto & [m_x],
\endaligned
$$
where $m_x$ denotes the maximal
ideal sheaf of $x$, is an isomomorphism. Fix a regular
point $x\in X$. Then we have an isomorphism $\phi_x\:\bar
J_{-1} @>>> \bar J_0$, by sending $[I]\in \bar J_{-1}$ to
$[I(x)]\in \bar J_0$. Under the composition $\psi:=\phi_x\circ
\Cal A^1$, the Cartier divisor $\Theta$ corresponds to the
Cartier divisor $[x]$ in $X$. 

Let $n\geq 3$ be an integer. The complete linear system
associated to $\Cal O_X(nx)$ gives rise to an embedding $X
\hookrightarrow \text{\bf P}^{n-1}$. If $H\subseteq \text{\bf
P}^{n-1}$ is a hyperplane intersecting $X$ at regular
points $y_1,\dots,y_n$, then $[y_1]+\dots+[y_n]$ is a Cartier
divisor on $X$ whose associated invertible sheaf is $\Cal
O_X(nx)$. Under $\psi$, the divisor 
$[y_1]+\dots+[y_n]$ corresponds to $\Theta_E$, where
$$
E=(\Cal O_X(y_1)\oplus\dots\oplus\Cal
O_X(y_n))\otimes\Cal O_X(-x).
$$
It follows now from Bertini's theorem that the theta sections of degree 
$n$ associated to completely decomposable vector bundles 
generate $H^0(\bar J_0, \Cal L^{\otimes n})$ for every $n\geq 0$. 
(In case $n\leq 2$ it is easy to check the latter statement directly.) 
Thus, for the case of curves of arithmetic genus 1, Question 9 
is answered in the affirmative. 

\vskip0.8cm

\parindent=0pt

\headline{\smallrm 10 \  \  \  \  \  \  \  \  \  \  \  
\  \  \  \  \  \  \  \  \  \  \  \  \  \  \  \  \  \  \  \  \  \  \  \  \  
\  \  \  \  \  \  \  \  \  \  \  \  \  \  \  \  \  \  \  \  \  \  \  \  \  
\  \  \  \  \  \  \  \  \  \  \  \  \  E. Esteves}

\smallit Acknowledgements. \smallrm 
This work is a natural continuation of part of the
author's Ph.D. thesis \cite{\smallrm 6}. 
The author would like to thank Prof. Kleiman
for introducing him to the theory of the compactified
Jacobian and for many useful comments. In addition, the author would like 
to thank Prof. Faltings for discovering a serious mistake in the 
first version of the present article. Finally, the author
would like to thank Waseda University, in particular
Prof. Kaji and his students, and Prof. Homma for their warm hospitality
during the period this work was carried out.

\headline{\smallrm Very ampleness for Theta on the compactified Jacobian 
\  \  \  \  \  \  \  \  \  \  \  \  \  \  \  \  \  \  \  \  \  
\  \  \  \  \  \  \  \  \  \  \  \  \  \  \  \  \  \  \  \  \  \  \  \  \  
\  \  \  \  \  \  \  \  \  \  \  \  \  \  \  \  \  11}

\enddocument